\documentclass[12pt]{article}
\usepackage{amsmath,amssymb}
\usepackage{rotating}
\usepackage{xypic}
\usepackage{amssymb}
\usepackage{amsmath,amscd}
\usepackage{pst-node,pstricks,multido,pst-plot,pst-text,pst-3d}%
\usepackage{graphicx}
\usepackage{amsfonts}
\newtheorem{example}{Example}[section]

\newtheorem{theorem}[example]{Theorem}

\newtheorem{corollary}[example]{Corollary}

\newtheorem{proposition}[example]{Proposition}

\def\S{{\mathfrak  S}}
\def\cal#1{{\mathfrak #1}}

\def\<{\langle}
\def\>{\rangle}

\def\C{{\mathbb C}}
\def\Z{{\mathbb Z}}

\def\ashuff#1#2#3{
\kern 1pt \vrule height#1 \overline{\vrule height#3 width 0pt
\hskip#2} \rule{.3pt}{#1}\overline{\vrule height#3 width 0pt
\hskip#2} \rule{.3pt}{#1} \kern 1pt }

\def\Det{{\rm Det}}

\def\Det{{\rm Det}}

\def\sign{\mbox{sign}}

\begin{document}

\title{Hyperdeterminantal computation for the Laughlin wave function}
\author{A. Boussicault, J.-G. Luque and C. Tollu}
 \maketitle
 \begin{abstract}
The decomposition of the Laughlin wave function in the Slater
orthogonal basis appears in the discussion on the second-quantized
form of the Laughlin states and is straightforwardly equivalent to
the decomposition of the even powers of the Vandermonde determinants
in the Schur basis. Such a computation is notoriously difficult and
the coefficients of the expansion have not yet been interpreted. In
our paper, we give an expression of these coefficients in terms of
hyperdeterminants of sparse tensors. We use this result to construct
an algorithm allowing to compute one coefficient of the development
without computing the others. Thanks to a program in {\tt C}, we
performed the calculation for the square of the Vandermonde up to an
alphabet of eleven lettres.
 \end{abstract}
 \section{Introduction}

When submitted to a magnetic field orthogonal to their motion,
electrons experience the Lorentz force which generates an asymmetric
distribution of the charge density in the conductor perpendicularly
to both  the line of sight path of the current and  the magnetic
field. The resulting voltage, called the Hall voltage, is
proportional to both the current and the magnetic flux density. To
extreme low temperature, in a strong magnetic field  and for a
two-dimensional electron system, the Hall conductance admits
quantized values which are integer or fractional multiples of
$e^2\over h$.
In the aim to explain this phenomenon, Laughlin \cite{Lau} proposed
quantum wave functions indexed by fractional fillings of the lowest
Landau level\footnote{Energy levels of a particule in a constant
uniform magnetic field \cite{LL}}. In the simplest cases
\cite{DGIL,Dun}, Fermi statistics require a fractional filling
$1\over 2k+1$ ($k$ being integer) and the corresponding Laughlin
wavefunction reads
\begin{equation}\begin{array}{rcl}
\Psi^{n,k}_{\rm
Laughlin}(z_1,\dots,z_n)&=&V(z_1,\dots,z_n)^{2k+1}\exp\{-\frac12\sum_{i=1}^N|z_i|^2\}\\
&=&V(z_1,\dots,z_n)^{2k}\Psi^0_{\rm
Laughlin}(z_1,\dots,z_n),\end{array}
\end{equation}
where $V(z_1,\dots,z_n)=\prod_{i<j}(z_i-z_j)$ is the Vandermonde
determinant.  Dunne \cite{Dun}\footnote{  Dunne discussed the
second-quantized form of the Laughlin states for the fractional
quantum Hall effect by decomposing the Laughlin wavefunctions into
the $n$-particle Slater basis and gives a general formula for the
expansion coefficients in terms of the characters of the symmetric
group.} and Di Francesco {\it et al.} \cite{DGIL} studied,
independently,  the expansion of the Laughlin wave function as a
linear combination of Slater wavefunctions for $n$ particules
\begin{equation}\label{SlaterDet}
\Psi_{\rm
Slater}^{\lambda}:=\frac1{\sqrt{n!\pi^n\prod_{i=1}^n\lambda_i!}}\exp\{-\frac12\sum_{i=1}^N|z_i|^2\}
\left|\begin{array}{cccc}
 z_1^{\lambda_1}&z_1^{\lambda_2}&\dots&z_1^{\lambda_n}\\
 z_2^{\lambda_1}&z_2^{\lambda_2}&\dots&z_2^{\lambda_n}\\
 \vdots&\vdots&&\vdots\\
 z_n^{\lambda_1}&z_n^{\lambda_2}&\dots&z_n^{\lambda_n}\end{array}\right|.
\end{equation}
It is easy to show that this problem is equivalent
 to the expansion of a power of the
discriminant in the Schur basis \cite{DGIL,Dun,KTW,STW}. Indeed, it
suffices to factorize the Slater wave function $\Psi_{\rm
Slater}^{\lambda}$ by the Schur function $S_\lambda$
\[
\Psi_{\rm
Slater}^{\lambda}=\frac1{\sqrt{n!\pi^n\prod_{i=1}^n\lambda_i!}}S_\lambda\Psi_{\rm
Laughlin}^0.
\]
A short time after the study of Di Francesco {\it et al.}, Sharf {\it
et al.} \cite{STW} proposed several algorithms to
compute this expansion. In particular, they performed it until $n=9$
for the square of the Vandermonde determinant and showed that a
conjecture (referred to as the admissibility condition) of \cite{DGIL} about the
characterization of the partitions having a non-null contribution in
the expansion fails for $n=8$. Note that King {\it et al.}
 showed \cite{KTW} that the conjecture becomes true if one considers
 the $q$-discriminant instead of the discriminant. In the same
 paper, they gave other methods for computing the expansion and
 perform it until $n=9$ in the case of the $q$-discriminant. In \cite{WSite}, the reader can found
  the expansion of $V^{2k}(z_1,\dots,z_n)$   until
 $n=10$ for $k=1$ and until $n=6$ for $k=2$.
 In the present paper, we give an expression of each coefficient as a hyperdeterminant
 (a natural generalization of the determinant for higher order
 tensors).
 As an application, we propose a new algorithm
  to compute each coefficient independently  from
 the others. The interest of such a result is twofold: First  the calculation
 can be distributed on several computers and the computation being essentially numerical, the algorithm can be implement
 in many programming languages. Second
 this method being based on the Laplace expansion of
hyperdeterminants, it allows us to write new recurrence formulae. 

\section{The Laughlin wavefunction and the admissibility conditions}
Di Francesco {\it et al.} \cite{DGIL} defined {\it admissible
partitions} as the partitions which {\it can appear} when one
expands $V(z_1,\dots,z_n)^{2k}$ on the Schur basis. That is the
partition arising as the dominant exponents when one expands
$V(z_1,\dots,z_n)^{2k+1}$ on the monomials without
simplifying. In other words, a partition $\lambda$ 
is admissible
if and only if $z^\lambda:=z^{\lambda_1}_1\dots z^{\lambda_n}_n$
appears with a nonvanishing coefficient in the expansion of
$\prod_{i<j}(z_i+z_j)^{2k+1}$,
\[
\prod_{ i<j}(z_i+z_j)^{2k+1}=\dots+ \alpha_\lambda z^\lambda +\dots\
..
\]

For a given pair of integers $n$ and $k$, the set of admissible
partitions is the interval for the dominance order ({\it i.e.}
$\lambda\geq \mu$ if and only if for each $0\leq i\leq n$,
$\lambda_1+\dots+\lambda_i\geq \mu_1+\dots+\mu_i$) whose upper bound
is $[2k(n-1),\dots,2k,0]$ and lower bound is
$[k(n-1),\dots,k(n-1)]$.

\noindent In \cite{DGIL}, Di Francesco {\it et al.} conjectured that
admissibility is a necessary and sufficient condition for
non-nullity of the coefficient $g_\lambda^{n,k}$. The first counter
example appears for $n=8$ and $k=1$ and was given by Scharf {\it et
al.} \cite{STW} who computed all the coefficients up to $n=9$, for
$k=1$.

\section{Hyperdeterminants}
\subsection{Definition and basics properties}
The birth of hyperdeterminants dates back to 1843, when Cayley gave
a lecture at the Cambridge Philosophical Society, about
functions that are reducible to sums of determinants. Actually,
Cayley used the same name of {\it
hyperdeterminant} to define several polynomials extending the notion of
determinants to higher order tensors. The polynomial which we use here, can be
considered as the simplest one because of its definition extending
in a natural way the expression of the determinant as an alternated
sum. Let ${\rm M}=(M_{i_1,\dots,i_p})_{1\leq i_1,\dots,i_p\leq n}$
be a tensor with $p$ indices, the {\it hyperdeterminant} of ${\rm M}$
is the alternated sum over $p$ copies of the symmetric group $\S_n$,
\begin{equation}\label{DefHypDet}
\Det({\rm
M}):=\frac1{n!}\sum_{\sigma_1,\dots,\sigma_p\in\S_n}\sign(\sigma_1\dots\sigma_p)\prod_{i=1}^nM_{\sigma_1(i)\dots\sigma_p(i)}.
\end{equation}
For example, if $p=4$ and $n=2$, \begin{multline*}\Det({\rm
M})=-M_{{2,1,1,1}}M_{{1,2,2,2}}+M_{{2,1,1,2}}M_{{1,2,2,1}}+M_{{2,1,2,1}}M_
{{1,2,1,2}}\\-M_{{2,1,2,2}}M_{{1,2,1,1}}+M_{{2,2,1,1}}M_{{1,1,2,2}}-M_{{
2,2,1,2}}M_{{1,1,2,1}}-M_{{2,2,2,1}}M_{{1,1,1,2}}+M_{{2,2,2,2}}M_{{1,1
,1,1}} \end{multline*} Straightforwardly,  $\Det$ is the zero
polynomial when $p$ is odd. Hence, we will suppose  that $p=2k$ is
even.

We will consider a special kind of hyperdeterminants: the {\it Hankel
hyperdeterminants}, whose entries depends only on the sums of the
indices,
\begin{equation}\label{DefHankel}
{\rm H}^f:=(f(i_1+\dots+i_{2k}))_{0\leq i_1,\dots,i_{2k}\leq n-1}.
\end{equation}
The Hankel hyperdeterminants appear in the literature in the work of
Lecat \cite{Lec2} (see also \cite{Lec1,Lec2}), but few properties
have been considered. More recently, one of the authors with
Jean-Yves Thibon \cite{LT1,LT2} and two of the authors with Hacene
Belbachir \cite{BBL}  investigated the links between these
polynomials and the Selberg integral and the Jack polynomials.

More generally, one defines a {\it shifted Hankel hyperdeterminant}
depending on $2k$ decreasing vectors
$\lambda^{(1)},\dots,\lambda^{(2k)}\in \Z^n$ as the hyperdeterminant
of the {\it shifted Hankel tensor}
\begin{equation}\label{DefSHankel}
{\rm
H}_{\lambda^{(1)},\dots,\lambda^{(2k)}}^f:=(f(\lambda_{i_1}^{(1)}+\dots+\lambda_{i_{2k}}^{(2k)}+i_1+\dots+i_{2k}))_{0\leq
i_1,\dots,i_{2k}\leq n-1}.
\end{equation}

\subsection{Minors of hypermatrices}
We will denote by  ${\rm M}\left[\begin{array}{c}
I_1\\\vdots\\I_{2k}\end{array}\right]$  the minor of a tensor ${\rm
M}=\left(M_{i_1,\dots,i_{2k}}\right)_{1\leq i_1,\dots,i_{2k}\leq n}$
obtained by choosing the elements indexed by the $2k$ increasing $m$-vectors $I_1,\dots, I_{2k}$, {\it i.e.,}
\[
{\rm M}\left[\begin{array}{c} I_1\\\vdots\\I_{2k}\end{array}\right]
:=\left(M_{j_{i_1}^{(1)},\dots,j_{i_{2k}}^{(2k)}}\right)_{1\leq
i_1,\dots,i_{2k}\leq m},\] if $I_1=(j_1^{(1)}\leq\dots\leq
j_{m}^{(1)})$, $\dots$, $I_{2k}=(j_1^{(2k)}\leq\dots\leq
j_{m}^{(2k)})$.

A shifted Hankel tensor is nothing but a minor of the infinite
Hankel tensor
\[
{\rm
H}^f_\infty:=\left(f(i_1+\dots+i_{2k}\right)_{-\infty<i_1,\dots,i_{2k}<\infty}.
\]

Hence, the property  {\it to be a shifted Hankel tensor} is closed
for the operation {\it extracting a minor}.

More generally,  consider the generic infinite tensor
\[
{\rm
M}_\infty:=\left(M_{i_1,\dots,i_{2k}}\right)_{-\infty<i_1,\dots,i_{2k}<\infty},
\]
and set for each $2k$-tuple of decreasing vectors
$\lambda^{(1)},\dots,\lambda^{(2k)}\in\Z^n$,
\[
{\rm M}_{\lambda^{(1)},\dots,\lambda^{(2k)}}=
\left(M_{n-\lambda^{(1)}_{i_1}+1+i_1,\dots,n-\lambda^{(2k)}_{i_{2k}}+1+i_{2k}}\right)_{1\leq
i_1,\dots,i_{2k}\leq n}.
\]
The tensor ${\rm M}_{\lambda^{(1)},\dots,\lambda^{(2k)}}$ is
obviously a minor of ${\rm M}_\infty$ and conversely, each minor
 ${\rm M}_\infty\left[\begin{array}{c}
I_1\\\vdots\\I_{2k}\end{array}\right]$ of ${\rm M}_\infty$ is equal
to some ${\rm M}_{\lambda^{(1)},\dots,\lambda^{(2k)}}$. Hence,
each minor of ${\rm M}_{\lambda^{(1)},\dots,\lambda^{(2k)}}$ is
again a minor of ${\rm M}_\infty$ and can be written in the form
${\rm M}_{\mu^{(1)},\dots,\mu^{(2k)}}$.

 More precisely,
one has the following property.

\begin{proposition}\label{PropMinors}(Compositions of minors)\\
Let $\lambda^{(1)},\dots,\lambda^{(2k)}\in\Z^n$ be $2k$ decreasing
vectors and $J_1,\dots,J_{2k}\subset\{1,\dots,n\}$ be $2k$ subsets
of $\{1,\dots,n\}$ with the same cardinality $m$, $0\leq m\leq n$. Then
the minor
\[
{\rm
M}_{\lambda^{(1)},\dots,\lambda^{(2k)}}\left[\begin{array}{c}\{1,\dots,n\}\setminus
J_1\\\vdots\\\{1,\dots,n\}\setminus J_{2k}
\end{array}\right]={\rm M}_{{\nu}^{(1)},\nu^{(2)},\dots,\nu^{(2k)}},
\]
where{\footnotesize
\begin{multline*}
\nu^{(p)}:=[\lambda^{(p)}_1+m,\dots,\lambda^{(p)}_{n-j_m}+m,\lambda^{(p)}_{n-j_m+2}+m-1,
\dots,\lambda^{(p)}_{n-j_{m-1}}+m-1,\lambda^{(p)}_{n-j_{m-1}+2}+m-2,\\\dots,\lambda^{(p)}_{n-j_1}+1,
\lambda^{(p)}_{n-j_1+2},\dots,\lambda_n^{(p)}]\end{multline*}} if
$J_p=\{j_1\leq\dots\leq j_m\}\subset\{1,\dots,n\}$.
\end{proposition}
{\bf Proof}
It suffices to understand the case of the vectors ({\it i.e.} the
tensors with only one indice). A straightforward induction on the
size of $J_1$ allows us to conclude. $\Box$

\subsection{A generalization of the Laplace expansion}
In the general case, there is no efficient algorithm for
computing an hyperdeterminant. Nevertheless, we will use a
generalization of the Laplace expansion for the hyperdeterminant due
to Zajaczkowski \cite{Zaj}\footnote{Armenante \cite{Arm} gave the
first generalization of the Laplace formula for cubic
hyperdeterminants whose first index is not alternating. Few years after
Zajaczkowski, Gegenbauer \cite{Geg} stated a Laplace formula for
general hyperdeterminants with alternating and not alternating indices.}.
\begin{theorem}\label{LaZajGeg}(Generalized Laplace) Zajaczkowski \cite{Zaj},
Gegenbauer \cite{Geg}\\ Consider a tensor ${\rm
M}=\left(M_{i_1,\dots,i_{2k}}\right)_{1\leq i_1,\dots,i_{2k}\leq
n}$, $0\leq m\leq n$ and $I_1=\{j_1^{(1)}\leq\dots\leq
j_m^{(1)}\}\subset\{1\dots n\}$. The hyperdeterminant of ${\rm M}$ can
be expanded as a alternated sum of $\left(n\atop m\right)^{2k-1}$
products of two minors,
\begin{equation}\label{Laplace}
\Det({\rm M})=\sum_{I_2,\dots,I_{2k}}\pm \Det\left({\rm
M}\left[\begin{array}{c}
I_1\\\vdots\\I_{2k}\end{array}\right]\right)\Det\left({\rm
M}\left[\begin{array}{c} \{1,\dots,n\}\backslash
I_1\\\vdots\\\{1,\dots,n\}\backslash I_{2k}\end{array}\right]\right)
\end{equation}
where the sum runs over the $m$-uplets $$I_2=[j_1^{(2)},\dots,
j_m^{(2)}], \,\dots, \,I_{2k}=[j_1^{(2k)},\dots, j_m^{(2k)}]
\in\{1,\dots,n\}^{m}$$ and $\pm$ denotes the sign of the product of
the permutations $\sigma_i$ bringing the indices of $I_i$ followed
by the indices of $\{1,\dots,n\}\backslash I_i$ into the original
order.
\end{theorem}

For general hypermatrices, the algorithm induced by this theorem is
not more efficient than the direct expansion but we will use it to
compute hyperdeterminants of sparse tensors.

\section{Computing the coefficients $g_\lambda^{n,k}$}
 \subsection{Hyperdeterminantal expression }
One can write some multiple integrals involving products of
determinants as hyperdeterminants.
\begin{proposition}\label{HeineG} (Generalized Heine identity)\\
Let $(f_j^{(i)})_{1\leq i\leq 2k\atop 1\leq j\leq n}$ be a family of
functions $\C\rightarrow\C$, and $\mu$ be any mesure on $\C$ such
that the integrals appearing in equality (\ref{dd2D}) are defined.
Then one has
\begin{multline}\label{dd2D}  {1\over n!}\int\dots\int \det(f_j^{(1)}(z_i))\dots
\det(f_j^{(2k)}(z_i))d\mu(z_1)\dots d\mu(z_n)=\\\Det\left(\int
f_{i_1}^{(1)}(z)\dots f_{i_{2k}}^{(2k)}(z)d\mu(z)\right)_{1\leq
i_1,\dots,i_{2k}\leq n}.
\end{multline}
\end{proposition}
{\bf Proof} Straightforward, expanding the left and right hand sides
of equality (\ref{dd2D}). $\Box$

In particular, if one applies Proposition \ref{HeineG} to the
product of a Schur function and a power of the discriminant, one
obtains a shifted Hankel hyperdeterminant whose entries are the
moments of the measure $\mu$.

\begin{corollary}Let $\mu$ be any mesure on $\C$ such
that the integrals appearing in equality (\ref{dd2D}) are defined.
One has {\footnotesize
\begin{equation}\label{ShiftHeine}
{1\over n!}\int\dots\int
s_\lambda(z_1,\dots,z_n)V(z_1,\dots,z_n)^{2k}d\mu(z_1)\dots
d\mu(z_n)=\Det(c_{\lambda_{n-i_1+1}+i_1+\dots+i_{2k}-2k})
\end{equation}}
where $c_n=\int z^nd\mu(z)$ denotes the $n$th moment of the measure
$\mu$.
\end{corollary}
{\bf Proof} It suffices to remark that
\[
s_\lambda(z_1,\dots,z_n)V(z_1,\dots,z_n)^{2k}=\det(z_i^{\lambda_{n-j+1}+j-1})\det(z_i^{j-1})^{2k-1},
\]
and to apply (\ref{dd2D}).$\Box$\\
 \noindent Let $\lambda^{(1)},\dots, \lambda^{(2k)}$ be $2k$ decreasing
vectors of $\Z^n$. One defines the tensor
\[
\Delta_{\lambda^{(1)},\dots,
\lambda^{(2k)}}:=\left(\delta_{\lambda^{(1)}_{n-i_1+1}+\dots+
\lambda^{(2k)}_{n-i_{2k}+1}+i_1+\dots+i_{2k},(2k-1)n+1}\right)_{1\leq
i_1,\dots,i_{2k}\leq n},
\]
and
\[
{\cal
D}_{\lambda^{(1)},\dots,\lambda^{(2k)}}:=\Det\left(\Delta_{\lambda^{(1)},\dots,
\lambda^{(2k)}}\right),
\]
its hyperdeterminant.

\noindent The following property gives an expression of the
coefficient $g_{\lambda}^{n,k}$ in terms of hyperdeterminants.

\begin{corollary} One has
\begin{equation}
g_{\lambda}^{n,k}=(-1)^{n(n-1)\over2} {\cal
D}_{\lambda,\underbrace{[0^n],\dots,[0^n]}_{2k+1\times}}
\end{equation}
\end{corollary}
{\bf Proof} This is a direct consequence of (\ref{ShiftHeine}) and
the definition of ${\cal D}_{\lambda^{(1)},\dots,\lambda^{(2k)}}$.
The complete discussion appears in \cite{BBL}.
 $\Box$

\subsection{Basic properties of the hyperdeterminants ${\cal
D}_{\lambda^{(1)},\dots,\lambda^{(2k)}}$} Let us list some
straightforward properties of such hyperdeterminants.
\begin{proposition}
Let $\sigma$ be any permutation of $\S_{2k}$, then
\[
{\cal D}_{\lambda^{(1)},\dots,\lambda^{(2k)}}={\cal
D}_{\lambda^{(\sigma(1))},\dots,\lambda^{(\sigma(2k))}}.
\]
\end{proposition}

\begin{proposition}
Let $m_1,m_2,\dots,m_{2k-1}\in\Z$  be $2k-1$ integers. One has,
\[
{\cal D}_{\lambda^{(1)},\dots,\lambda^{(2k)}}={\cal
D}_{[\lambda^{(1)}_1+m_1,\dots,\lambda^{(1)}_n+m_1],\dots,
[\lambda^{(2k)}_1+m_{2k},\dots,\lambda^{(2k)}_n+m_{2k}]},
\]
where $m_{2k}=-m_1-\dots-m_{2k-1}$.
\end{proposition}

\begin{proposition}
If ${\cal D}_{\lambda^{(1)},\dots,\lambda^{(2k)}}\neq 0$ then
\[
\sum \lambda_i^{(j)}=(k-1)n(n-1).
\]
\end{proposition}
{\bf Proof} If $ {\cal D}_{\lambda^{(1)},\dots,\lambda^{(2k)}}\neq
0$ then there exist $\sigma_1,\dots,\sigma_{2k}\in\S_n$ such that
\[
\prod_{i=1}^n
\delta_{\lambda^{(1)}_{n-\sigma_1(i)+1}+\dots+\lambda^{(2k)}_{n-\sigma_{2k}(i)+1}+\sigma_1(i)+\dots+\sigma_{2k}(i),(2k-1)n+1}
\neq 0.
\]

This implies
\[
\sum_i\lambda^{(1)}_{n-\sigma_1(i)+1}+\dots+\lambda^{(2k)}_{n-\sigma_{2k}(i)+1}+\sigma_1(i)+\dots+\sigma_{2k}(i)=n((2k-1)n+1).
\]
But the left hand side is nothing but
$\sum_{i,j}\lambda_j^{(i)}+kn(n+1).$ The result follows. $\Box$

\subsection{Minors of the matrices $\Delta_{\lambda^{(1)},\dots,\lambda^{(2k)}}$}
 Consider the sets defined by
 \[
  \Gamma_{k,n}:=\{\Delta_{\lambda^{(1)},\dots,\lambda^{(2k)}}|\lambda^{(1)},\dots,\lambda^{(2k)}\mbox{
  are decreasing vectors of }\Z^n\}.
 \]

\begin{proposition}
Let $\lambda^{(1)},\dots,\lambda^{(2k)}\in\Z^n$ be $2k$ decreasing
vectors and $J_1,\dots,J_{2k}\subset\{1,\dots,n\}$ be $2k$ subsets
of $\{1,\dots,n\}$ with the same cardinality $m$, $0\leq m\leq n$. Hence
the minor
\[
\Delta_{\lambda^{(1)},\dots,\lambda^{(n)}}\left[\begin{array}{c}\{1,\dots,n\}\setminus
J_1\\\vdots\\\{1,\dots,n\}\setminus J_{2k}
\end{array}\right]
\]
belongs to $\Gamma_{k,n-m}$.
\end{proposition}
{\bf Proof} From Proposition \ref{PropMinors}, one obtains
\begin{multline*}
\Delta_{\lambda^{(1)},\dots,\lambda^{(2k)}}\left[\begin{array}{c}\{1,\dots,n\}\setminus
J_1\\\vdots\\\{1,\dots,n\}\setminus J_{2k}
\end{array}\right]=\\
\left(\delta_{\nu_{n-i_1+1}^{(1)}+\dots+\nu_{n-i_{2k}+1}^{(2k)}+i_1+\dots+i_{2k},(2k-1)(n-m+1)+1}\right)_{1\leq
i_1,\dots,i_{2k}\leq n-m}
\end{multline*}
 where
 {\footnotesize
\begin{multline*}
\nu^{(p)}:=[\lambda^{(p)}_1+m,\dots,\lambda^{(p)}_{n-j_m}+m,\lambda^{(p)}_{n-j_m+2}+m-1,
\dots,\lambda^{(p)}_{n-j_{m-1}}+m-1,\lambda^{(p)}_{n-j_{m-1}+2}+m-2,\\\dots,\lambda^{(p)}_{n-j_1}+1,
\lambda^{(p)}_{n-j_1+2},\dots,\lambda_n^{(p)}]\end{multline*}} if
$J_p=\{j_1\leq\dots\leq j_m\}\subset\{1,\dots,n\}$. Furthermore,
{\footnotesize
\begin{multline*}
\delta_{\nu_{n-i_1+1}^{(1)}+\dots+\nu_{n-i_{2k}+1}^{(2k)}+i_1+\dots+i_{2k},(2k-1)(m+1)+1}=
\delta_{{\nu'}_{n-i_1+1}^{(1)}+\dots+\nu_{n-i_{2k}+1}^{(2k)}+i_1+\dots+i_{2k},(2k-1)(n-m)+1},
\end{multline*}}
where ${\nu'}^{(1)}$ is the decreasing sequence
\begin{equation}\label{Dminor}{\nu'}^{(1)}:=[{\nu_1}^{(1)}-m(2k-1),\dots,{\nu_n}^{(1)}-m(2k-1)].
\end{equation} Hence,
\[
\Delta_{\lambda^{(1)},\dots,\lambda^{(2k)}}\left[\begin{array}{c}\{1,\dots,n\}\setminus
J_1\\\vdots\\\{1,\dots,n\}\setminus J_{2k}
\end{array}\right]=\Delta_{{\nu'}^{(1)},\nu^{(2)},\dots,\nu^{(2k)}}\in\Gamma_{k,n-m}.
\]
This completes the proof.
 $\Box$

\subsection{A recursive formula for ${\cal
D}_{\lambda^{(1)},\dots,\lambda^{(2k)}}$}

As a consequence of the preceding sections, one has

\begin{corollary}
 Let $1\leq p\leq m$, one has
\begin{equation}\label{recD}
{\cal
D}_{\lambda^{(1)},\dots,\lambda^{(2k)}}=\sum_{I}(-1)^{i_1+\dots+i_{2k}}{\cal
D}_{\mu^{(1)}_I,\dots,\mu^{(2k)}_I},
\end{equation}
where  the sum is over the $2k$-tuples,
$I=[i_1,i_2,\dots,i_{2k}]\in\{1,\dots,n\}^{2k}$ verifying
\[
\lambda^{(1)}_{n-i_1+1}+\dots+\lambda^{(2k)}_{n-i_{2k}+1}+i_{1}+\dots+i_{2k}=(2k-1)n+1,
\]
and the decreasing vectors $\mu^{(1)}_I,\dots, \mu^{(2k)}_I$ are
defined by
\begin{multline}\label{mu}
\\
\mu_I^{(1)}=[\lambda_1^{(1)}-2(k-1),\dots,\lambda_{n-i_1}^{(1)}-2(k-1),\lambda_{n-i_1+2}^{(1)}-2(k-1)-1,
\dots,\lambda_{n}^{(1)}-2(k-1)-1]\\
\mu_I^{(2)}=[\lambda_1^{(2)}+1,\dots,\lambda_{n-i_2}^{(2)}+1,\lambda_{n-i_2+2}^{(2)},
\dots,\lambda_{n}^{(2)}]\\
\vdots\\
\mu_I^{(2k)}=[\lambda_1^{(2k)}+1,\dots,\lambda_{n-i_{2k}}^{(2k)}+1,\lambda_{n-i_{2k}+2}^{(2k)},
\dots,\lambda_{n}^{(2k)}].
\end{multline}
\end{corollary}
{\bf Proof} Setting $I_1=\{i_1\}$ in Theorem \ref{LaZajGeg}, from
the definition of ${\cal D  }_{\lambda^{(1)},\dots,\lambda^{(2k)}}$,
one gets,
$${\cal
D}_{\lambda^{(1)},\dots,\lambda^{(2k)}}=\sum_{I}(-1)^{i_1+\dots+i_{2k}}
\Det\left(\Delta_{\lambda^{(1)},\dots,\lambda^{(k)}}
\left[\begin{array}{c}\{1,\dots,n\}\setminus
i_1\\\vdots\\\{1,\dots,n\}\setminus i_{2k}\end{array}\right]
\right)$$ where the sum is over the $2k$-tuples,
$I=[i_1,i_2,\dots,i_{2k}]\in\{1,\dots,n\}^{2k}$ verifying
\[
\lambda^{(1)}_{n-i_1+1}+\dots+\lambda^{(2k)}_{n-i_{2k}+1}+i_{1}+\dots+i_{2k}=(2k-1)n+1.
\]
Furthermore, one has
$$
\Delta_{\lambda^{(1)},\dots,\lambda^{(2k)}}
\left[\begin{array}{c}\{1,\dots,n\}\setminus
i_1\\\vdots\\\{1,\dots,n\}\setminus
i_{2k}\end{array}\right]=\Delta_{\mu^{(1)},\dots,\mu^{(2k)}}
$$
where the partitions $\mu^{(i)}$ are defined by (\ref{mu}). The
result follows.
 $\Box$

\begin{example}\rm
Suppose that we want to compute ${\cal D}_{[211][100][100][000]}$. That
is to compute the hyperdeterminant of
\[
\Delta_{[211][100][100][000]}=\tiny \begin{array}{c|ccccccccc|}
_{i_3,i4}\setminus^{i1,i2}&11&\tiny 12&\tiny 13&\tiny 21&\tiny
22&\tiny 23&\tiny 31&\tiny 32&\tiny 33\\
\hline
 \tiny 11& & & & & & & & & \\
 \tiny 12& & & & & &1& &1& \\
 \tiny 13& & &1& & & &1& & \\
 \tiny 21& & & & & &1& &1& \\
 \tiny 22& & &1& & & &1& & \\
 \tiny 23& & & & &1& & & & \\
 \tiny 31& & & & &1& & & & \\
 \tiny 32& &1& &1& & & & & \\
 \tiny 33&1& & & & & & & &
\end{array}.
\]
If one sets $i_1=1$, the only indices $(i_1,i_2,i_3,i_4)$ such that
the corresponding entries of $\Delta_{[211][100][100][000]}$ do not
vanish are $(1,3,2,2), (1,2,3,2), (1,1,3,3)$  and $(1,3,1,3)$.
Furthermore
\[
\Delta_{[211][100][100][000]}\left[\begin{array}{c}
\{2,3\}\\\{1,2\}\\\{1,3\}\\\{1,3\}\end{array}\right]=\tiny
\begin{array}{c|cccc|} _{i_3,i4}\setminus^{i1,i2}&
   11&12&21&22\\
\hline
11& & & & \\
12& & &1& \\
21& &1& &  \\
22& & & &
\end{array}=\Delta_{[0 -1][00][20][10]}.
\]

On the same way, one has
\[
\Delta_{[211][100][100][000]}\left[\begin{array}{c}
\{2,3\}\\\{1,3\}\\\{1,2\}\\\{1,3\}\end{array}\right]=\Delta_{[0
-1][20][00][10]},
\]
\[
\Delta_{[211][100][100][000]}\left[\begin{array}{c}
\{2,3\}\\\{2,3\}\\\{1,2\}\\\{1,2\}\end{array}\right]=\Delta_{[0
-1][21][00][00]}
\]
and
\[
\Delta_{[211][100][100][000]}\left[\begin{array}{c}
\{2,3\}\\\{1,2\}\\\{2,3\}\\\{1,2\}\end{array}\right]=\Delta_{[0
-1][00][21][00]}.
\]
Hence,
\[
{\cal D}_{[211][100][100][000]}={\cal D}_{[0 -1][00][20][10]}+{\cal
D}_{[0 -1][20][00][10]}+{\cal D}_{[0 -1][21][00][00]}+{\cal D}_{[0
-1][00][21][00]}.\] A straightforward computation gives
\[{\cal D}_{[0
-1][20][00][10]}={\cal D}_{[0 -1][00][20][10]}=1\] and
\[
{\cal D}_{[0 -1][21][00][00]}={\cal D}_{[0 -1][00][21][00]}=2.
\]
from what it follows that ${\cal D}_{[211][100][100][000]}=6$.
\end{example}

\begin{example}\rm Here one illustrates the fact that the recurrence
(\ref{recD}) provides an algorithm to compute the coefficient
$g^{k,n}_\lambda$ . Suppose that one wants to compute the
coefficient of $s_{411}$ in the square of the Vandermonde
determinant for an alphabet of size $3$. One needs to compute the
value of ${\cal D}_{[4,1,1],[0,0,0],[0,0,0],[0,0,0]}$. Applying the
Laplace expansion, one finds  that this can be written as a sum
involving $27$ hyperdeterminants

\begin{multline*}
{\cal D}_{[4,1,1],[0,0,0],[0,0,0],[0,0,0]}=\alpha_{3111}{\cal
D}_{[2,-1],[0,0],[0,0],[0,0]}+\alpha_{3112}{\cal
D}_{[2,-1],[0,0],[0,0],[1,0]}+\\\dots+\alpha_{3333}{\cal
D}_{[2,-1],[1,1],[1,1],[1,1]}. \end{multline*} But for only three of
them the coefficient $\alpha_{I}$ does not vanish
\begin{multline*}
{\cal D}_{[4,1,1],[0,0,0],[0,0,0],[0,0,0]}=\alpha_{3211}{\cal
D}_{[2,-1],[1,0],[0,0],[0,0]}+\alpha_{3112}{\cal
D}_{[2,-1],[0,0],[0,0],[1,0]}\\+\alpha_{3121}{\cal
D}_{[2,-1],[0,0],[1,0],[0,0]}. \end{multline*} One has
$\alpha_{3112}=\alpha_{3121}=\alpha_{3211}=-1$ and for reason of
symmetry $${\cal D}_{[2,-1],[1,0],[0,0],[0,0]}={\cal
D}_{[2,-1],[0,0],[0,0],[1,0]}={\cal D}_{[2,-1],[0,0],[1,0],[0,0]}.$$
It remains to compute ${\cal D}_{[2,-1],[1,0],[0,0],[0,0]}$. Using
again the Laplace expansion, one finds that this can be written as
the sum of $8$ hyperdeterminants, of which only one gives a nonvanishing $\alpha_I$,
$${\cal D}_{[2,-1],[1,0],[0,0],[0,0]}=\alpha_{2111}{\cal D}_{[0],[0],[0],[0]}=-1$$
Hence, $g_{411}^{1,3}=3$.
\end{example}

\subsection{Factorisation formul\ae}

\begin{proposition}\label{FactD}
Let $\lambda^{(1)},\dots,\lambda^{2k}$ such that there exists an
integer $0<m<n$ verifying
\[
\lambda^{(1)}_{1}+\dots+\lambda^{(1)}_m+\lambda_{n-m}^{(2)}+\dots+\lambda_{n}^{(2)}
+\dots+\lambda_{n-m}^{(2k)}+\dots+\lambda_{n}^{(2k)}=(k-1)m(m-1)
\]
then ${\cal D}_{\lambda^{(1)},\dots,\lambda^{(2k)}}$ factorizes as

\[
{\cal D}_{\lambda^{(1)},\dots,\lambda^{(2k)}}=\pm{\cal
D}_{\mu^{(1)},\dots,\mu^{(2k)}}{\cal
D}_{\nu^{(1)},\dots,\nu^{(2k)}},
\]
where
\begin{multline*}
\mu^{(1)}:=[\lambda^{(1)}_1-2(k-1)m,\dots,\lambda^{(1)}_m-2(k-1)m],
\
\nu^{(1)}:=[\lambda^{(1)}_{m+1},\dots,\lambda^{(1)}_n]\\
\mu^{(2)}:=[\lambda^{(2)}_{n-m+1},\dots,\lambda^{(2)}_n], \
\nu^{(2)}:=[\lambda^{(2)}_{1},\dots,\lambda^{(2)}_{n-m}]\\
\vdots\\
\mu^{(2k)}:=[\lambda^{(2k)}_{n-m+1},\dots,\lambda^{(2k)}_n], \
\nu^{(2k)}:=[\lambda^{(2k)}_{1},\dots,\lambda^{(2k)}_{n-m}]
\end{multline*}
\end{proposition}
{\bf Proof} It is a direct consequence of the generalized Laplace
expansion.
 $\Box$
\begin{corollary}\label{Factg}
Let $\lambda$ be such that it exists an integer $0<m<n$ verifying
\[
\lambda_{1}+\dots+\lambda_m=km(m-1)
\]
then
\[
g_{\lambda}^{n,k}= g_{\mu}^{n-m,k}g_{\nu}^{m,k},
\]
where
\[\mu:=[\lambda_1-2k(m-1),\dots,\lambda_m-2k(m-1)],
\mbox{ and } \nu:=[\lambda_{m+1},\dots,\lambda_n].\]
\end{corollary}
{\bf Proof} It is a direct consequence of Proposition
\ref{FactD}.$\Box$

Note that Corollary \ref{Factg} can also be obtained as a
straightforward consequence of the factorization
\[
\Delta(x_1,\dots,x_n)=\Delta(x_1,\dots,x_m)\prod_{i=1}^m\prod_{j=m+1}^n(x_i-x_j)\Delta(x_{m+1},\dots,x_n).
\]
\begin{example}
 To calculate the coefficient $g_{77420}^{1,5}$, one may compute the
 hyperdeterminant
$${\cal D}_{[77420],[0000],[0000],[0000]}.$$ From Proposition
\ref{FactD}, it factorizes as
$$\pm{\cal D}_{[420],[000],[000],[000]}{\cal D}_{[33],[00],[00],[00]} ..$$
Hence,
$$g_{77420}^{1,5}=g_{420}^{1,3}g_{33}^{1,2}. $$
\end{example}

\section{Results}
The rules explained in the previous sections enables to write an
algorithm computing the coefficients $g_\lambda^{k,n}$. The
calculations being completely numerical, they can be implemented in a
programming language such as {\tt C} which allows us to optimize
runtime and  memory management. A program written in {\tt C} can be
downloadeed from \cite{Site}. All calculations have been
performed on a personal computer\footnote{Intel Pentium processor
1.86Ghz, 1Go Ram.}, with the only exception of the case $k=1$ and $n=11$,
for which a 8-processors cluster with 32 Go Ram was used. In  the most
general case, computing a hyperdeterminant using the generalized
Laplace theorem is possible only for very small dimensions. Here, as
we consider only very sparse tensors, the computation can be
achieved for reasonably large alphabets.\\ Table \ref{T1}  contains
the list of the cases which have been computed with this program.
\begin{table}[h]\
\begin{equation*}
\begin{array}{|c|c|}
\hline k&n\ max\\
\hline 1&\mbox{up to }11\\
2&\mbox{up to }7\\
3&\mbox{up to }6\\
4&\mbox{up to }5\\
5&\mbox{up to }5\\\hline
\end{array}
\end{equation*}
\caption{\label{T1} List of the case for which the computation have
been performed for all admissible partitions.}
\end{table}
The results can be downloaded from \cite{Site}. As expected, there are
fewer nonvanishing partitions than admissible partitions. Tables
\ref{T2} and \ref{T3} contain respectively the number of admissible
partitions and the number of vanishing admissible partitions.
\begin{table}[h]
{\footnotesize$$\
\begin{array}{|c|cccccccccc|}
\hline
\begin{array}{l}\end{array}&
n=2&3&4&5&6&7&8&9&10&11\\
\hline k=1& 2&5&16&59&247&1111&5302&28376&135670&716542\\
k=2&3&13&76&521&3996&32923&&&&\\
k=3&4&25&213&2131&23729&&&&&\\
k=4&5&41&459&6033&88055&&&&&\\
k=5&6&61&846&13771&&&&&&\\\hline
\end{array}
$$
} \caption{\label{T2} Number of admissible partitions}
\end{table}
\begin{table}[h]
 {\footnotesize$$\
\begin{array}{|c|cccccccccc|}
\hline &
n=2&3&4&5&6&7&8&9&10&11\\
\hline k=1& 0&0&0&0&0&0&8&66&389&1671\\
k=2&0&0&0&0&6&46&?&?&?&?\\
k=3&0&0&0&2&14&?&?&?&?&?\\
k=4&0&0&0&16&?&?&?&?&?&?\\
k=5&0&0&0&0&?&?&?&?&?&?\\\hline
\end{array}
$$
} \caption{\label{T3} Number of vanishing admissible partitions}
\end{table}
 \section{Conclusion}
We have described an algorithm which computes each coefficient
appearing in the expansion of the Laughlin wave functions in the
Slater basis without computing the others, which allows to distribute
easily the computation. This algorithm is based on an interpretation
of each coefficient as an hyperdeterminant. This approach being
completely numerical our algorithm can be implemented in various
languages (such as {\tt C}). 
The principal limitation of our method is that the generalization
to the $q$-deformation is not easy. In particular, one has to
construct an analogue of the (multi)-antisymmetrizer. One
possible approach would consist in searching for the latter operator in the double
affine Hecke algebra. Indeed, in previous articles, two of the
authors gave $q$-deformations \cite{BL,qdef} which can  be written
as symmetric Macdonald functions
 indexed by rectangular or staircase partitions for some
  specializations of the parameters (which
made us think that the  Hecke algebra may play a r\^ole). We have not identified
the operator yet.\\ The method can also be
adapted to write the powers of the discriminant in the monomial
basis. In this case, one has to compute the hypedeterminant with
non alternating indices of a sparse tensor using the more general
version of the Gegenbauer-Laplace expansion theorem \cite{Geg}.
Nevertheless, the tensor considered are bigger (with an odd number
of indices). Furthermore, several others methods exist to perform
this computation (see {\it e.g.} \cite{STW}) and we do not know whether
ours is very efficient in this case.

It is also worth noting that Physicists use another and more efficient method to
carry out these calculations. They proceed by diagonalization of the
unphysical model Hamiltonian for which the power of the Vandermonde
is the exact ground state (see {\it e.g.} \cite{RGJ,RJ,RRC}).
The drawback of that algorithm is that one cannot obtain one
coefficient without computing the others. Another advantage of our
method is that it is based on a combinatorial description of some
hyperdeterminants (after recoding them, one only uses the vectors
which index them). Giving new relations, this  can be used to
understand the very difficult problem of the  characterization of
the partitions which have a nonvanishing contribution. One can follow
two tracks to solve this problem. The first one consists in
understanding the combinatorics of these hyperdeterminants. The
second, more algebraic and geometric, consists in characterizing the
varieties defined by the vanishing of a hyperdeterminant.

Finally, as the powers of the Vandermonde are special cases  of the
Read-Rezayi states \cite{RR}, one can naturally ask the question of
the generalization of our method to  other cases.

\noindent{\bf Acknowledgment}\\
Two of the authors (J.-G.L. and A.B) are grateful to Th. Jolicoeur
for useful discussions about the fractional quantum Hall effect.
J.-G.L. is grateful to C. Toke for discussions on the Read-Rezayi
states.

 \end{document}